\magnification=1200
\overfullrule=0pt
{\centerline {{\bf Some research perspectives in nonlinear functional
analysis}}\par
\bigskip
\bigskip
{\centerline {BIAGIO RICCERI}}
\bigskip
\bigskip
\bigskip
\bigskip
The object of this lecture is to propose a series of conjectures and problems
in different fields of analysis. They have been formulated with the aim
of introducing some innovative methods in the study of classical topics, as
open mappings, fixed points, critical points, global minima, control theory.
\par
\smallskip
We start recalling the following definition.
Let $(E,\|\cdot\|)$ be a real normed space. A non-empty set $A\subset
E$ is said to be {\it antiproximinal with respect to $\|\cdot\|$} if, for
every $x\in E\setminus A$ and every $y\in A$, one has
$$\|x-y\|>\inf_{z\in A}\|x-z\|\ .$$\par
\bigskip
CONJECTURE 1. -  There exists a non-complete real normed space $E$
with the following property: for every non-empty convex set $A\subset E$
which is antiproximinal with respect to each norm on $E$, the interior
of the closure of $A$ is non-empty.\par
\bigskip
The main reason for the study of Conjecture 1 is to give a contribution
to open mapping theory in the setting of non-complete normed spaces.
Actually, making use of Theorem 4 of [8], one can prove the following
result.\par
\bigskip
THEOREM 1. - {\it Let $X, E$ be two real vector spaces, $C$ a non-empty convex
subset of $X$, $F$ a multifunction from $C$ onto $E$, with non-empty values
and convex graph.\par
Then, for every non-empty convex set $A\subseteq C$  which is open with
respect to the relativization to $C$ of the strongest vector topology on $X$,
the set $F(A)$ is antiproximinal with respect to each norm on $E$.}\par
\bigskip
We now present\par
\bigskip
CONJECTURE 2. -  Let $E$ be a real Banach space, and let $J:E\to {\bf R}$
be a continuously G\^ateaux differentiable functional. Assume that there
are $r>0$ and $x_{0}, x_{1}\in E$, with $\|x_{0}-x_{1}\|>r$, such that
$$\inf_{\|x-x_{0}\|=r}J(x)\geq \max\{J(x_{0}),J(x_{1})\}\ .$$
Put
$$c=\inf_{u\in {\cal A}}\sup_{t\in [0,1]}J(u(t))$$
where ${\cal A}$ denotes the set of all continuous functions $u:[0,1]\to E$
such that $u(0)=x_{0}$, $u(1)=x_{1}$.\par
Then, for each $\epsilon>0$, there exist an interval $I\subset [0,1]$,
a $u\in {\cal A}$ and a continuous function $\varphi:I\times E\to {\bf R}
$, with $\varphi(t,\cdot)$ $\epsilon$-Lipschitzian in $E$ for all $t\in I$,
 such that
$$\sup_{t\in I}|J(u(t))-c|<\epsilon$$
and the set
$$\{(t,y)\in I\times E : J'(u(t))(y)=\varphi(t,y)\}$$
is disconnected. \par
\bigskip
CONJECTURE 3. - Let $(H,\langle\cdot,\cdot\rangle)$ be a real Hilbert
space, $X\subset H$ a non-empty compact convex set, $f:X\to X$ a continuous
function different from the identity, $\epsilon$ a positive real number
small enough. Denote by $\Lambda_{\epsilon}$ the set of all continuous
function $\varphi:X\times H\to {\bf R}$ such that, for each $x\in X$,
$\varphi(x,\cdot)$ is Lipschitzian in $H$, with Lipschitz constant less
than or equal to $\epsilon$. Consider $\Lambda_{\epsilon}$ equipped with the
relativization of the strongest vector topology on the space ${\bf R}^
{X\times H}$.\par
Then, the set
$$\{(\varphi,x,y)\in \Lambda_{\epsilon}\times X\times H:
\langle f(x)-x,y\rangle = \varphi(x,y)\}$$
is disconnected.\par
\bigskip
Conjectures 2 and 3 have been formulated with the aim of finding drastically
novel proofs of the mountain pass lemma and the Brower fixed point theorem,
respectively. Actually, if they were true, then we could get the mentioned
proofs by means of the following results from [15] (see also [12]):\par
\bigskip
THEOREM 2 ([15], Theorem 19). - {\it Let $X$ be a connected topological space,
 $E$ a real Banach space (with topological dual space $E^{*}$),
  $A$ a continuous operator from $X$ into $E^{*}$, $\varphi$
a real function on $X\times E$ such that,
for each $x\in X$, $\varphi(x,\cdot)$ is Lipschitzian in $E$, with 
Lipschitz constant $L(x)\geq 0$. Further, assume that the set
$$\{(x,y)\in X\times E :  A(x)(y) =
\varphi(x,y)\}$$
is disconnected.\par
Then, there exists some $v\in X$ such that
$\parallel A(v)\parallel_{E^{*}}\leq L(v).$}\par
\bigskip
THEOREM 3 ([15], Theorem 21). - {\it Let $X$ be a connected and locally
connected
topological space, $E$ a real Banach space, $A : X\rightarrow
E^{*}$ a continuous operator with closed range. For
each $\epsilon>0$, denote by $\Lambda_{\epsilon}$
the set of all continuous functions $\varphi:X\times E
\rightarrow {\bf R}$ such that, for each $x\in X$, $\varphi(x,\cdot)$
is Lipschitzian in $E$, with Lipschitz constant less than or equal
to $\epsilon$.
Consider $\Lambda_{\epsilon}$ equipped with the relativization of the
strongest vector topology on the space ${\bf R}^{X\times E}$, and
assume that the set
$$\{(\varphi,x,y)\in \Lambda_{\epsilon}\times X\times E :
 A(x)(y) =\varphi(x,y)\}$$
is disconnected.\par
Then, $A$ vanishes at some point of $X$.}\par
\bigskip
We now come to a problem in control theory.
Let $a$ be a positive real number and let $F$ be a given
multifunction from $[0,a]\times {\bf R}^{n}$ into ${\bf R}^{n}$.
We denote by ${\cal S}_{F}$ the set of all Carath\'eodory solutions
of the problem $x'\in F(t,x), x(0)=0$ in $[0,a]$. That is to
say
$${\cal S}_{F}=\{u\in AC([0,a],{\bf R}^{n}):
u'(t)\in F(t,u(t))\hskip 10pt \hbox {\rm a.e. in}\hskip 4pt
 [0,a],\hskip 3pt u(0)=0\}\ .$$
For each $t\in [0,a]$, put
$${\cal A}_{F}(t)=\{u(t) : u\in {\cal S}_{F}\}.$$
In other words, ${\cal A}_{F}(t)$ denotes the attainable
 set at time $t$. Also,
put
$$V_{F}=\bigcup_{t\in [0,a]}{\cal A}_{F}(t).$$
Finally, set
$$C_{F}=\{x\in {\bf R}^{n} : \{t\in [0,a] : x\in {\cal A}_{F}(t)\}
\hskip 7pt  \hbox {\rm is connected}\}.$$\par
\bigskip
PROBLEM 1. - Find conditions under which the set $C_{F}$ is non-empty and
open.
\bigskip
The study of the above problem would be interesting in view of the following
result, where "dim" stands for covering dimension:\par
\bigskip
THEOREM 4 ([11], Theorem 9). - {\it Assume that $F$
has non-empty compact convex values and bounded 
range. Moreover, assume that $F(\cdot,x)$ is measurable
for each $x\in {\bf R}^{n}$ and that $F(t,\cdot)$ is
upper semicontinuous for a.e. $t\in [0,a]$.\par
 Then, for every non-empty 
connected set
$X\subseteq V_{F}\cap C_{F}$ which is open in \hbox {\rm aff}$(X)$
and different from $\{0\}$, 
   one has the following alternative:\par
\noindent
 either
$$X\subseteq {\cal A}_{F}(a)$$ or
 $$\dim({\cal A}_{F}(t)\cap X)\geq \dim(X)-1$$
for some $t\in ]0,a[$.}\par
\bigskip
Now, we are going to present a problem about a unusual way of finding global
minima of functionals in Banach spaces. First, we recall the following result
from [17] (see also [10], [13], [14], [18], [19]):
\bigskip
THEOREM 5 ([17], Theorem 2.1) . - {\it
Let  $(T,{\cal F},\mu)$ be 
non-atomic measure space, with $\mu(T)<+\infty$,
  $E$ a real Banach space, and $f:E\to {\bf R}$
 a bounded below Borel functional such that
$$\sup_{x\in E}{{f(x)}\over {\|x\|^{\gamma}+1}}<+\infty$$
for some $\gamma\in ]0,1[$.\par
Then, for every $p\geq 1$ and every closed hyperplane $V$ of
 $L^{p}(T,E)$, one has
$$\inf_{u\in V}\int_{T}f(u(t))d\mu=\inf_{u\in L^{p}(T,E)}
\int_{T}f(u(t))d\mu.$$}
\bigskip
PROBLEM 2. - Let $X$ be an infinite-dimensional
real Banach space, and let $J:X\to {\bf R}$ be a
bounded below functional satisfying
$$\sup_{u\in X}{{J(u)}\over {\|u\|^{\gamma}+1}}<+\infty
\eqno{(*)}$$
for some $\gamma\in ]0,1[$.\par
Find conditions under which there exists a closed hyperplane $V$ of $X$
such that the restriction of $J$ to $V$ has a local minimum.\par
\bigskip
The motivation for the study of Problem 2 is as follows. Assume that we
wish to minimize a bounded below Borel functional $f$
on a real Banach space $E$
satisfying the condition
$$\sup_{x\in E}{{f(x)}\over {\|x\|^{\gamma}+1}}<+\infty$$
for some $\gamma\in ]0,1[$.\par
For each $u\in L^{1}([0,1],E)$, put
$$J(u)=\int_{0}^{1}f(u(t))dt\ .$$
So, $J$ is bounded below and satisfies $(*)$ with $X=L^{1}([0,1],E)$.
Assume that there is some closed hyperplane $V$ of $L^{1}([0,1],E)$ such
that the restriction of $J$ to $V$ has a local minimum, say $u_{0}$. By a
result of Giner ([3]) $u_{0}$ is actually a global minimum of the restriction of
$J$ to $V$. On the other hand, by Theorem 5, we have
$$\inf_{u\in V}J(u)=\inf_{u\in L^{1}([0,1],E)}J(u)$$
and so $u_{0}$ is a global minimum of $J$ in $L^{1}([0,1],E)$. This easily
implies that $f$ has a global minimum in $E$.\par
\bigskip
Using a major tool adopted in the proof of Theorem 5, O. Naselli got the
following wonderful characterization:\par
\bigskip
THEOREM 6 ([7], Theorem 1). - 
  {\it
Let  $(T,{\cal F},\mu)$ be a $\sigma$-finite
non-atomic complete measure space, $X$ a real topological vector
space, $\Phi$  a linear homeomorphism from $X$ 
onto $L^{1}(T)$, and  $f:T \times {\bf R} \to {\bf R} $ a
Carath\'eodory function.\par
Then, the following are equivalent:\par
\noindent
{\hbox {\rm (a)}}\hskip 5pt The set
$$\{u \in X: f(t, \Phi (u)(t)) = 0 \phantom{P} a.e.\hskip 3pt in
\hskip 3pt T \} $$ 
 is arcwise connected and intersects each closed hyperplane
of $X$.\par
\noindent
{\hbox {\rm (b)}}\hskip 5pt The function $$t \to \inf
 \{ \vert x \vert : x \in {\bf R},\hskip 3pt
 f(t,x) = 0 \} $$ 
belongs to
$L^{1}(T)$ and, for almost every $t \in T$, one has 
$$\sup \{ x \in {\bf R} : f(t,x) = 0 \} = + \infty $$
and 
$$\inf \{ x \in {\bf R} : f(t,x) = 0 \} = - \infty .$$ }\par
\bigskip
Surprisingly enough, the following problem seems to be open:\par
\bigskip
PROBLEM 3. Does Theorem 6 hold replacing $L^{1}(T)$ by $L^{p}(T)$ with
$p>1$ ?\par
\bigskip
The next problems concern a very particular function space introduced in
[9].
Let $m,n$ be two positive integers. Denote by $V({\bf R}^n)$
the space of all functions $u\in C^{\infty}({\bf R}^n)$ such that, for each
bounded subset $\Omega\subset {\bf R}^n$, one has
$$\sup_{\alpha\in {\bf N}_{0}^{n}}\sup_{x\in \Omega}|D^{\alpha}u(x)|
<+\infty,$$
where $D^{\alpha}u=\partial^{\alpha_{1}+...+\alpha_{n}}u/
\partial x_{1}^{\alpha_{1}}...\partial x_{n}^{\alpha_{n}}$,
$\alpha=(\alpha_{1},...,\alpha_{n})$ and ${\bf N}_{0}=
{\bf N}\cup \{0\}.$\par
\bigskip
PROBLEM 4. - Let $f:{\bf R}\to {\bf R}$ be a function such that, for each
$u\in V({\bf R}^n)$, the composite
function $x\to f(u(x))$ belongs to $V({\bf R}^n)$.
\par
Then, must $f$ necessarily be of the form $f(t)=at+b$ ?\par
\bigskip
PROBLEM 5. -
For each $\alpha\in {\bf N}_{0}^{n}$, with
$|\alpha|=\alpha_{1}+...+\alpha_{n}\leq m$, let
$a_{\alpha}\in {\bf R}$ be given.
 Let $P:V({\bf R}^n)\to V({\bf R}^n)$ be the differential operator defined
by putting
$$P(u)=\sum_{|\alpha|\leq m}a_{\alpha}D^{\alpha}u$$
for all $u\in V({\bf R}^n)$.\par
 Find necessary and sufficient conditions in order that 
$$P(V({\bf R}^n))=V({\bf R}^n).$$
Up to date, the only (very partial) answer to Problem 5 is provided by
the following
\bigskip
THEOREM 7 ([9], Theorem 4). - {\it Let $a, b\in {\bf R}\setminus \{0\}$ and
$h, k\in {\bf N}$. For each $u\in V({\bf R}^2)$, put
$$P(u)=a{{\partial^{h} u}\over {\partial x^{h}}}+
b{{\partial^{k} u}\over {\partial y^{k}}}.$$
Then, one has
$$P(V({\bf R}^2))=V({\bf R}^2)$$
if and only if $|a|\neq |b|$.} \par
\bigskip
The final problems we want to present come from specific applications to
nonlinear boundary value problems of the following result obtained in
[20]. \par
\bigskip
THEOREM 8 ([20], Theorem 2.5). - {\it Let $X$ be a non-empty sequentially
weakly closed set in a reflexive
 real Banach space, and let $\Phi, \Psi
:X\to ]-\infty,+\infty]$ be two sequentially weakly lower semicontinuous
functionals.
 Assume also that $\Psi$ is
(strongly) continuous. Denote by $I$ the set of all real numbers $\rho>
\inf_{X}\Psi$ such that $\Psi^{-1}(]-\infty,\rho[)$ is bounded and
intersects the domain of $\Phi$. Assume that $I\neq \emptyset$.
 For each $\rho\in I$, put 
$$\varphi(\rho)=\inf_{x\in \Psi^{-1}(]-\infty,\rho[)}
{{\Phi(x)-\inf_{\overline {(\Psi^{-1}(]-\infty,\rho[))}_{w}}\Phi}
\over {\rho-\Psi(x)}}\ ,$$
where $\overline {(\Psi^{-1}(]-\infty,\rho[))}_{w}$ is the closure
of $\Psi^{-1}(]-\infty,\rho[)$ in the relative weak topology of $X$.
 Furthermore, set
$$\gamma=\liminf_{\rho\to (\sup I)^{-}}\varphi(\rho)$$
and
$$\delta=\liminf_{\rho\to 
\left ( {\inf\limits _{X}\Psi}\right ) ^{+}}\varphi(\rho)\ .$$
Then, the following conclusions hold:\par
\smallskip
\noindent
$(a)$\hskip 10pt For each $\rho\in I$ and each $\mu>
\varphi(\rho)$, the restriction of the
functional $\Phi+\mu\Psi$ to $\Psi^{-1}(]-\infty,\rho[)$ has a global
minimum.\par
\smallskip
\noindent
$(b)$\hskip 10pt If $\gamma<+\infty$, then, for each $\mu>
\gamma$, the following alternative holds: either the restriction of
$\Phi+\mu\Psi$ to $\Psi^{-1}(]-\infty,\sup I[)$ has a global minimum,
 or there exists a sequence $\{x_{n}\}$
of local minima of $\Phi+\mu\Psi$ such that $\Psi(x_{n})<\sup I$ for
all $n\in {\bf N}$ and
$\lim_{n\to \infty}\Psi(x_{n})=\sup I$.\par
\smallskip
\noindent
$(c)$\hskip 10pt If $\delta<+\infty$, then, for each $\mu>
\delta$, 
 there exists a
sequence $\{x_{n}\}$ of local minima of
$\Phi+\mu\Psi$, with $\lim_{n\to \infty}\Psi(x_{n})=\inf_{X}\Psi$,
 which weakly converges to a global minimum of
$\Psi$.}\par
\bigskip
From now on, $\Omega$ is an open bounded subset of ${\bf R}^n$, with smooth
boundary, and (for $p>1$)
$W^{1,p}(\Omega)$, $W^{1,p}_{0}(\Omega)$ are the usual Sobolev
spaces, with norms
$$\|u\|=\left ( \int_{\Omega}|\nabla u(x)|^{p}dx + \int_{\Omega}|u(x)|^{p}dx
\right ) ^{1\over p}$$
and
$$\|u\|=\left ( \int_{\Omega}|\nabla u(x)|^{p}dx\right ) ^{1\over p}$$
respectively.\par
Let $p>1$, and  let $:\Omega\times {\bf R}\to {\bf R}$ be
 a Carath\'eodory function.\par
\smallskip
Recall that a weak solution of the Dirichlet problem
$$\cases {-\hbox {\rm div}(|\nabla u|^{p-2}\nabla u)
=f(x,u)
 & in
$\Omega$\cr & \cr u_{|\partial \Omega}=0\cr}$$
 is any $u\in W^{1,p}_{0}(\Omega)$ such that
 $$\int_{\Omega}|\nabla u(x)|^{p-2}\nabla u(x)\nabla v(x)dx
-\int_{\Omega}f(x,u(x))v(x)dx=0 \eqno{(**})$$
for all $v\in W^{1,p}_{0}(\Omega)$. While, a weak solution of the Neumann
problem
$$\cases {-\hbox {\rm div}(|\nabla u|^{p-2}\nabla u)=
f(x,u)
 & in
$\Omega$\cr & \cr {{\partial u}\over {\partial \nu}}=0 & on
$\partial \Omega$\cr} $$
 $\nu$ being the outer unit normal to $\partial \Omega$, is
 is any $u\in W^{1,p}(\Omega)$ satisfying identity $(**)$
for all $v\in W^{1,p}(\Omega)$. \par
\bigskip
Let us recall the following classical result by Ambrosetti and Rabinowitz:\par
\bigskip
THEOREM 9 ([1], Theorem 3.10). - {\it Assume that:\par
\noindent
$(1)$\hskip 5pt there are two positive constants $a, q$, with
$q<{{n+2}\over {n-2}}$ if $n\geq 3$, such that
$$|f(x,\xi)|\leq a(1+|\xi|^{q})$$
for all $(x,\xi)\in \Omega\times {\bf R}$;\par
\noindent
$(2)$\hskip 5pt there are constants $r\geq 0$ and $c>2$ such that
$$0<c\int_{0}^{\xi}f(x,t)dt\leq \xi f(x,\xi)$$
for all $(x,\xi)\in \Omega\times {\bf R}$ with $|\xi|\geq r$;\par
\noindent
$(3)$\hskip 5pt one has
$$\lim_{\xi\to 0}{{f(x,\xi)}\over {\xi}}=0$$   
uniformly with respect to $x$.\par
Then, the problem
$$\cases {-\Delta u=f(x,u)
 & in
$\Omega$\cr & \cr u_{|\partial \Omega}=0.\cr}$$
has a non zero weak solution.}\par
\bigskip
What can be said if, in Theorem 9, condition $(3)$ is removed at all ? Using
Theorem 8 (part $(a)$), we got the following result:\par
\bigskip
THEOREM 10 ([21, Theorem 4). - {\it Assume that conditions $(1)$ and $(2)$
hold.\par
Then, for each $\rho>0$ and each $\mu$
satisfying
$$\mu>\inf_{u\in B_{\rho}}{{\sup_{v\in B_{\rho}}\int_{\Omega}
\left ( \int_{0}^{v(x)}f(x,\xi)d\xi\right )dx-
\int_{\Omega}
\left ( \int_{0}^{u(x)}f(x,\xi)d\xi\right )dx}\over
{\rho-\int_{\Omega}|\nabla u(x)|^{2}dx}}\ , \eqno{(***)}$$
where
$$B_{\rho}=\left \{ u\in W^{1,2}_{0}(\Omega):
\int_{\Omega}|\nabla u(x)|^{2}dx<\rho\right \}\  ,$$
the problem
$$\cases {-\Delta u={{1}\over {2\mu}} f(x,u)
 & in
$\Omega$\cr & \cr u_{|\partial \Omega}=0\cr}$$
has at least two weak solutions one of which lies in $B_{\rho}$.}\par
\bigskip
The following problem naturally arises in connection with Theorems 9 and
10.\par
\bigskip
PROBLEM 6. -  Under conditions $(1)$ and $(2)$, is there some $\rho>0$ such
that the infimum appearing
in $(***)$ is less than ${{1}\over {2}}$ ? \par
\bigskip
Clearly, if the answer to this problem was positive, then Theorem 10 would be
a proper improvement of Theorem 9.\par
\bigskip
Again applying Theorem 8 (part $(a)$), in [23], we obtained the following
bifurcation theorem:\par
\bigskip
THEOREM 11. - {\it Let $f, g:\Omega\times {\bf R}\to {\bf R}$
be two Carath\'eodory functions. Assume that:\par
\noindent
$(i)$\hskip 10pt there is $s>1$ such that
$$\limsup_{\xi\to 0^+}
{{\sup_{x\in \Omega}|f(x,\xi)|}\over
{\xi^{s}}}<+\infty\ ;$$
\noindent
$(ii)$\hskip 10pt there is $q\in ]0,1[$ such that\par
$$\limsup_{\xi\to 0^+}{{\sup_{x\in \Omega}
|g(x,\xi)|}\over {\xi^q}}<+\infty\ ;$$
\noindent
$(iii)$\hskip 10pt there are a non-empty open set $D\subseteq \Omega$
and a set $B\subseteq D$ of positive measure such that
$$ \limsup_{\xi\to 0^{+}}{{\inf_{x\in B}
\int_{0}^{\xi}g(x,t)dt}\over {\xi^2}}=+\infty\ ,\hskip 3pt
 \liminf_{\xi\to 0^{+}}{{\inf_{x\in D}
\int_{0}^{\xi}g(x,t)dt}\over {\xi^2}}>-\infty\ .$$
Then, for some $\lambda^{*}>0$ and for
each $\lambda\in ]0,\lambda^{*}[$, 
the problem  
$$\cases {-\Delta u=
f(x,u)+\lambda g(x,u)
 & in
$\Omega$\cr & \cr u_{|\partial \Omega}=0\ , \cr} \eqno{\left ( P_{\lambda}
\right ) }$$
admits a non-zero, non-negative weak solution $u_{\lambda}\in
C^{1}(\overline {\Omega})$. Moreover, one has
$$\limsup_{\lambda\to 0^+}{{\|u_{\lambda}\|_{C^{1}(\overline {\Omega})}}\over
{\lambda^{q\over 1-q}}}<+\infty$$
and the function
$$\lambda\to {{1}\over {2}}\int_{\Omega}|\nabla u_{\lambda}(x)|^{2}dx-
\int_{\Omega}\left ( \int_{0}^{u_{\lambda}(x)}f(x,\xi)d\xi\right ) dx
-\lambda\int_{\Omega}\left ( \int_{0}^{u_{\lambda}(x)}g(x,\xi)d\xi
\right ) dx$$
 is negative
and decreasing in $]0,\lambda^{*}[$.
If, in addition,
 $f, g$ are continuous in $\Omega\times [0,+\infty[$
and
$$\liminf_{\xi\to 0^{+}}{{\inf_{x\in \Omega}g(x,\xi)}\over
{\xi |\log\xi|^{2}}}>-\infty\ ,$$
then $u_{\lambda}$ is positive in $\Omega$.}\par
\bigskip
In view of [2], where problem $\left ( P_{\lambda}\right )$ is studied for
particular nonlinearities, we point out the following problem:\par
\bigskip
PROBLEM 7. - Under the assumptions of Theorem 11, does
problem $\left ( P_{\lambda}\right )$ admit a non-zero,
non-negative, {\it minimal} solution
for each $\lambda>0$ small enough ?    \par
\bigskip
In the two final theorems, $\lambda$ denotes a function in
$L^{\infty}(\Omega)$, with $\hbox {\rm ess inf}_{\Omega}\lambda>0$. They have
been established in [22] as applications of parts $(b)$ and $(c)$ of Theorem
8. Their non-smooth versions have been obtained in [4] and [6].\par
\bigskip
THEOREM 12 ([22], Theorem 3). - {\it Assume $p>n$.
Let $f:{\bf R}\to {\bf R}$ be a continuous function,
and $\{a_{k}\}$, $\{b_{k}\}$ two sequences
in ${\bf R}^+$ satisfying
$$a_{k}<b_{k}\hskip 5pt \forall k\in {\bf N},\hskip 5pt
\lim_{k\to \infty}b_{k}=+\infty,\hskip 5pt 
\lim_{k\to\infty}{{a_{k}}\over {b_{k}}}=0\ ,$$
$$\max\left \{ \sup_{\xi\in [a_{k},b_{k}]}\int_{a_{k}}^{\xi}
f(t)dt, \sup_{\xi\in [-b_{k},-a_{k}]}\int_{-a_{k}}^{\xi}
f(t)dt\right \}\leq 0\hskip 5pt \forall k\in {\bf N}$$
and
$$\limsup_{|\xi|\to +\infty}{{\int_{0}^{\xi}f(t)dt}\over
{|\xi|^{p}}}=+\infty\ .$$
Then, for every $\alpha, \beta\in L^{1}(\Omega)$, with $\min\{\alpha(x),
\beta(x)\}\geq 0$ a.e. in $\Omega$ and $\alpha\neq 0$,
and for every continuous function $g:{\bf R}\to {\bf R}$ satisfying
$$\sup_{\xi\in {\bf R}}\int_{0}^{\xi}g(t)dt\leq 0$$
and
$$\liminf_{|\xi|\to +\infty}{{\int_{0}^{\xi}g(t)}\over {|\xi|^{p}}}>
-\infty\ ,$$
the problem
$$\cases {-\hbox {\rm div}(|\nabla u|^{p-2}\nabla u)+\lambda(x) |u|^{p-2}u=
\alpha(x)f(u)+\beta(x)g(u) & in
$\Omega$\cr & \cr {{\partial u}\over {\partial \nu}}=0 & on
$\partial \Omega$\cr} $$
admits an unbounded sequence of weak solutions in $W^{1,p}(\Omega)$.}
\bigskip
THEOREM 13 ([22], Theorem 4). - {\it Assume $p>n$.
Let $f:{\bf R}\to {\bf R}$ be a continuous function,
and $\{a_{k}\}$, $\{b_{k}\}$ two sequences
in ${\bf R}^+$ satisfying
$$a_{k}<b_{k}\hskip 5pt \forall k\in {\bf N},\hskip 5pt
\lim_{k\to \infty}b_{k}=0,\hskip 5pt 
\lim_{k\to\infty}{{a_{k}}\over {b_{k}}}=0\ ,$$
$$\max\left \{ \sup_{\xi\in [a_{k},b_{k}]}\int_{a_{k}}^{\xi}
f(t)dt, \sup_{\xi\in [-b_{k},-a_{k}]}\int_{-a_{k}}^{\xi}
f(t)dt\right \}\leq 0\hskip 5pt \forall k\in {\bf N}$$
and
$$\limsup_{\xi\to 0}{{\int_{0}^{\xi}f(t)dt}\over
{|\xi|^{p}}}=+\infty\ .$$
Then, for every $\alpha, \beta\in L^{1}(\Omega)$, with $\min\{\alpha(x),
\beta(x)\}\geq 0$ a.e. in $\Omega$ and $\alpha\neq 0$,
and for every continuous function $g:{\bf R}\to {\bf R}$ satisfying
$$\sup_{\xi\in {\bf R}}\int_{0}^{\xi}g(t)dt\leq 0$$
and
$$\liminf_{\xi\to 0}{{\int_{0}^{\xi}g(t)}\over {|\xi|^{p}}}>
-\infty\ ,$$
the problem
$$\cases {-\hbox {\rm div}(|\nabla u|^{p-2}\nabla u)+\lambda(x) |u|^{p-2}u=
\alpha(x)f(u)+\beta(x)g(u) & in
$\Omega$\cr & \cr {{\partial u}\over {\partial \nu}}=0 & on
$\partial \Omega$\cr} $$
admits a sequence of non zero weak solutions which strongly
converges to $0$ in $W^{1,p}(\Omega)$.}\par
\bigskip
A major problem about Theorems 12 and 13 is as follows:\par
\bigskip
PROBLEM 8. - In Theorems 12 and 13, when $g=0$, are the conclusions still
valid without the assumption
$$\lim_{k\to\infty}{{a_{k}}\over {b_{k}}}=0\hskip 3pt ?$$\par
\bigskip
A partial answer to this problem (for Theorem 13)
has recently been provided by G. Anello and
G. Cordaro in [3].\par
\vfill \eject
{\centerline {{\bf References}}}\par
\bigskip
\bigskip
\noindent
[1]\hskip 5pt A. AMBROSETTI and P. H. RABINOWITZ, {\it Dual variational
methods in critical point theory and applications}, J. Funct. Anal.,
{\bf 14} (1973), 349-381.\par
\smallskip
\noindent
[2]\hskip 5pt A. AMBROSETTI, H. BREZIS and G. CERAMI, {\it Combined
effects of concave and convex nonlinearities in some elliptic problems},
J. Funct. Anal., {\bf 122} (1994), 519-543.\par
\smallskip
\noindent
[3]\hskip 5pt G. ANELLO and G. CORDARO. {\it Infinitely many positive
solutions for the Neumann problem involving the $p$-Laplacian}, preprint.
\par
\smallskip
\noindent
[4]\hskip 5pt P. CANDITO, {\it Infinitely many solutions to the Neumann
problem for elliptic equations involving the $p$-Laplacian and with
discontinuous nonlinearities}, Proc. Math. Soc. Edinburgh, to appear.\par
\smallskip
\noindent
[5]\hskip 5pt E. GINER, {\it Minima sous contrainte, de fonctionnelles
int\'egrales}, C. R. Acad. Sci. Paris, S\'erie I, {\bf 321} (1995),
429-431.\par
\smallskip
\noindent
[6]\hskip 5pt S. A. MARANO and D. MOTREANU, {\it Infinitely many critical
points of non-differentiable functions and applications to a Neumann type
problem involving the $p$-Laplacian}, J. Differential Equations, to appear.
\par
\smallskip
\noindent
[7]\hskip 5pt O. NASELLI, {\it On the solution set of an equation
of the type $f(t,\Phi(u)(t))=0$}, Set-Valued Anal., {\bf 4} (1996),
399-405.\par
\smallskip
\noindent
[8]\hskip 5pt B. RICCERI, {\it Images of open sets under certain
multifunctions}, Rend. Accad. Naz. Sci. XL, {\bf 10} (1986), 33-37.\par
\smallskip
\noindent
[9]\hskip 5pt B. RICCERI, {\it On the well-posedness of the Cauchy problem
for a class of linear partial differential equations of infinite order in
Banach spaces}, J. Fac. Sci. Univ. Tokyo, Sec. IA, {\bf 38}
 (1991), 623-640.\par
\smallskip
\noindent
[10]\hskip 5pt B. RICCERI, {\it A variational property of integral
functionals on $L^{p}$-spaces of vector-valued functions}, 
C. R. Acad. Sci. Paris, S\'erie I, {\bf 318} (1994), 337-342.\par
\smallskip
\noindent
[11]\hskip 5pt B. RICCERI, {\it Applications of a theorem concerning sets
with connected sections}, Topol. Methods Nonlinear Anal.,
{\bf 5} (1995), 237-248.\par
\smallskip
\noindent
[12]\hskip 5pt B. RICCERI, {\it Existence of zeros via disconnectedness},
J. Convex Anal., {\bf 2} (1995), 287-290.\par
\smallskip
\noindent
[13]\hskip 5pt B. RICCERI, {\it A variational property of integral functionals
and related conjectures}, Banach Center Publ., {\bf 35} (1996), 237-242.\par
\smallskip
\noindent
[14]\hskip 5pt B. RICCERI, {\it On the integrable selections of certain
multifunctions}, Set-Valued Anal., {\bf 4} (1996), 91-99.\par
\smallskip
\noindent
[15]\hskip 5pt B. RICCERI, {\it Recent uses of connectedness in functional
analysis}, RIMS, Kyoto, Surikai\-
sekikenkyusho-Kokyuroku, {\bf 939} (1996), 11-22.\par
\smallskip
\noindent
[16]\hskip 5pt B. RICCERI, {\it On some motivated conjectures and problems},
 Matematiche, {\bf 51} (1996), 369-373.\par
\smallskip
\noindent
[17]\hskip 5pt B. RICCERI, {\it More on a variational property of
integral functionals}, J. Optim. Theory Appl.,
{\bf 94} (1997), 757-763.\par
\smallskip
\noindent
[18]\hskip 5pt B. RICCERI, {\it On a topological minimax theorem and
its applications}, in ``Minimax theory and applications'', B. Ricceri
and S. Simons eds., 191-216, Kluwer Academic Publishers, 1998.\par
\smallskip
\noindent 
[19]\hskip 5pt B. RICCERI, {\it Further considerations on a variational
property of integral functionals}, J. Optim. Theory Appl., {\bf 106}
(2000), 677-681.\par
\smallskip
\noindent
[20]\hskip 5pt B. RICCERI, {\it A general variational principle and
some of its applications}, J. Comput. Appl. Math., {\bf 113}
(2000), 401-410.\par
\smallskip
\noindent
[21]\hskip 5pt B. RICCERI, {\it On a classical existence
theorem for nonlinear elliptic equations}, in ``Experimental,
constructive and nonlinear analysis'', M. Th\'era ed.,
275-278, CMS Conf. Proc. {\bf 27}, Canad. Math. Soc., 2000.\par
\smallskip
\noindent
[22]\hskip 5pt B. RICCERI, {\it Infinitely many solutions of the Neumann
problem for elliptic equations involving the p-Laplacian},
Bull. London Math. Soc., {\bf 33} (2001), 331-340.\par
\smallskip
\noindent
[23]\hskip 5pt B. RICCERI, {\it A bifurcation theorem for nonlinear
elliptic equations}, preprint. \par
\smallskip
\noindent
\bigskip
\bigskip
Department of Mathematics\par
University of Catania\par
Viale A. Doria 6\par
95125 Catania, Italy
\bye